\documentclass[12pt]{article}
\usepackage{amsmath, amsfonts, amssymb, amsthm}
\usepackage{epsfig}
\usepackage{multirow}

\let\BFseries\bfseries\def\bfseries{\BFseries\mathversion{bold}} 

\theoremstyle{plain}
\newtheorem{thm}{Theorem}
\newtheorem{lem}[thm]{Lemma}

\newtheorem{cor}[thm]{Corollary}
\theoremstyle{definition}

\newtheorem{rem}[thm]{Remark}

\renewenvironment{proof}[1][] {\smallskip \noindent {\bf Proof#1:} }{\hspace*{\fill}$\square$\medskip\par}
\def\P{{\bf {\mathbb{P}}}}
\newcommand{\pr}[1]{\P\left(#1\right)}

\def\E{\mathbb{E}}

\newcommand{\deq}{\stackrel{d}{=}}

\begin{document}
\title{Wait-and-see strategies in polling models}
\author{Frank Aurzada, Sergej Beck, and Michael Scheutzow}
\date{\today}
\maketitle
\begin{abstract} We consider a general polling model with $N$ stations. The stations are served exhaustively and in cyclic order. Once a station queue falls empty, the server does not immediately switch to the next station. Rather, it waits at the station for the possible arrival of new work (``wait-and-see'') and, in the case of this happening, it restarts service in an exhaustive fashion. The total time the server waits idly is set to be a fixed, deterministic parameter for each station. Switchover times and service times are allowed to follow some general distribution, respectively. In some cases, which can be characterised, this strategy yields strictly lower average queueing delay than for the exhaustive strategy, which corresponds to setting the ``wait-and-see credit'' equal to zero for all stations. This extends results of Pek\"{o}z (Probability in the Engineering and Informational Sciences 13 (1999)) and of Boxma et al.\ (Annals of Operations Research 112 (2002)). Furthermore, we give a lower bound for the delay for {\it all} strategies that allow the server to wait at the stations even though no work is present.
\end{abstract}

\bigskip


\bigskip

\noindent{\slshape\bfseries Authors' address:}\\  Technische Universit\"{a}t Berlin, Institut f\"{u}r  Mathematik, Sekr.\ MA 7-5, Stra\ss e des 17.\ Juni 136, 10623 Berlin, Germany.\\ aurzada@math.tu-berlin.de, stahanovez@gmail.com, ms@math.tu-berlin.de,\\ phone: +493031424219, fax: +493031421695.

\section{Introduction and main results}
\subsection{Introduction}
In this work, we consider a polling model in the sense of \cite{takagi}. In a polling model, one server serves several queues, called stations. The classical service procedures are the (a) exhaustive, (b) gated, and (c) limited strategies, where the server serves each station (a) until no more work is waiting at the respective station, (b) until all the work is served that was awaiting the server upon its arrival at the station, or (c) until the server has finished at most a predescribed number of jobs. The server then turns its attention to the next station. A possible (deterministic or random) idle time between the different stations, called switchover time, accounts for things like reloading or refueling.

Recently, a few papers (most importantly \cite{boxma2002} and \cite{pekoz1999}, also see \cite{boxma2002b} and \cite{yechiali}) consider strategies where the server does not immediately switch from one station to the next if the queue there is empty. Rather, it possibly waits at the station for a while for the potential arrival of new messages. This is particularly useful if (i) the switchover times are random with sufficiently large variances and (ii) if the server is not likely to find much work at the other stations, that is, if the traffic intensity of the current station is much larger than those of the other stations.

The strategy proposed and analysed in this paper is also of this type: Each station $i$ is given a fixed wait-and-see credit $T_i\geq 0$. Once the server arrives at station $i$, it will work there whenever messages are waiting, but it will also wait (and see) at the station for a total time of $T_i$. Once the credit is used up and no more messages are waiting, it will switch to the next station. This strategy was considered by Pek\"{o}z \cite{pekoz1999} for the case of a completely symmetric system, where all arrival rates, service times, switchover times, and the $T_i$ are identical for all stations.

The main contributions of this paper are
\begin{itemize}
 \item to extend Pek\"{o}z' results to the general (that is, not necessarily symmetric) polling model and to show in particular that the asymmetry induces some new effects previously not observed,
 \item to show that our strategy can be adjusted to provide lower delay than with the exhaustive strategy in several cases (which can be characterised and which also appear for deterministic switchover times),
 \item to analyse the case of a polling model with two stations in detail and compare our strategy to the one proposed by Boxma et al.\ \cite{boxma2002}, and
 \item to prove a lower bound for the delay for {\it all} strategies that allow the server to wait at a station even though no work may be present.
\end{itemize}

As mentioned above, introducing a wait-and-see credit is particularly useful if the server is not likely to find much traffic at the other stations. This is because changing the station means to stay idle for a switchover time rather than resuming work at the current station within a short time. Surprisingly, we will see that this effect is largely independent of the length of the switchover times.

So far, the advantage of additional idle times -- as we apply them here -- was ascribed to the random switchover times. The new observation is that using non-zero idle times is also particularly useful if the system is asymmetric, that is, one of the stations experiences much more traffic than the others, even though the switchover times may be even deterministic. This is an aspect that could not be observed in \cite{pekoz1999}, and even though being intuitive, we can quantify this effect precisely.

In our polling model, the stations are served in cyclic order. We mention that the performance of all strategies can yet be improved by altering the order in which the server serves the different stations. For example, star polling can be applied if one of the stations experiences significantly more traffic than all other stations (see e.g.\ \cite{boxma1991,olsen}).

This paper is structured as follows. In Section~\ref{sec:model}, we describe the model in detail and introduce the relevant parameters. Section~\ref{sec:main} contains a summary of our main results. We review related work in Section~\ref{sec:relwork}. The proofs for the main results are given in Section~\ref{sec:general} for the general case and in Section~\ref{sec:two} for the refined results for polling models with only two stations. In Section~\ref{sec:lower}, we prove a lower bound for the delay for {\it all} strategies that allow the server to wait at a station even though no work may be present. We highlight some possible further improvements and lines of future reseach in Section~\ref{sec:outlook}.

The motivation for this work comes from a real world application: In so-called Ethernet Passive Optical Networks (EPONs, see \cite{kramer,lam}), a service provider is connected to various end users via an optical fibre cable. Different optical wavelength channels may be available on the cable for the communication, but each wavelength channel can be operated only either upstream (messages are sent from end users to the service provider) or downstream at a given time. Switching from upstream to downstream operation or vice versa incurs an idle time (switchover time). Therefore, each channel of an EPON can be regarded as a polling model.

\subsection{The model}\label{sec:model}
We consider a polling model with $N\geq 1$ stations and one server which serves the stations in cyclic order. The stations are numbered $i=1,\ldots, N$; because of the cyclic order, when we talk of the stations, we set $N+1\triangleq 1$.

Each station $i$ has its own queue which is fed by a Poisson arrival process whose intensity is denoted by $\lambda_i$. Each arriving message has a random length (also called service time). The mean and second moment of the message length distribution are denoted by $b_i$ and $b_i^{(2)}$, respectively, and are assumed to be finite.

The behaviour of the server can be described as follows. The server arrives at station $i$ and starts serving (FCFS) all waiting messages and newly arriving messages until the queue is empty. This is typically called exhaustive service in the context of polling models. However, once the station is empty or if the server finds an empty station upon its arrival, the server does not immediately switch to the next station; it rather turns idle for some time in order to wait for potentially newly arriving messages (``wait-and-see''). As soon as new messages arrive, it starts serving them immediately and in an exhaustive fashion. Once finished, it again turns idle and waits for new messages to arrive, and so on.

The main feature of our model is that the server is set to wait idly for new messages for a \emph{total} time $T_i$, where $T_i\geq 0$ is a fixed parameter of the system, called wait-and-see credit. This total time can be spent altogether in one single period, for example, if there are no messages waiting at the station upon the server's arrival and no messages arrive even until time $T_i$ after the server's arrival at the station; or it can be spent in different periods -- interleaved by different busy periods. Note that, since $T_i$ is fixed, the server may not use any information about the current queue status at other stations nor about the future of the arrival process at any station.

After the server has spent a total waiting time of $T_i$ at station $i$, it starts the switchover to station $i+1$. Hereby, it first spends a possibly random idle time, called switchover time, where it does not serve any messages neither at station $i$ nor at station $i+1$. The random switchover time from station $i$ to station $i+1$ is assumed to have finite mean $r_i\geq 0$ and finite second moment $r_i^{(2)}$. We will consider both non-deterministic and deterministic switchover times (in the latter case $r_i^{(2)}=r_i^2$).

\begin{figure}
\includegraphics[scale=0.88]{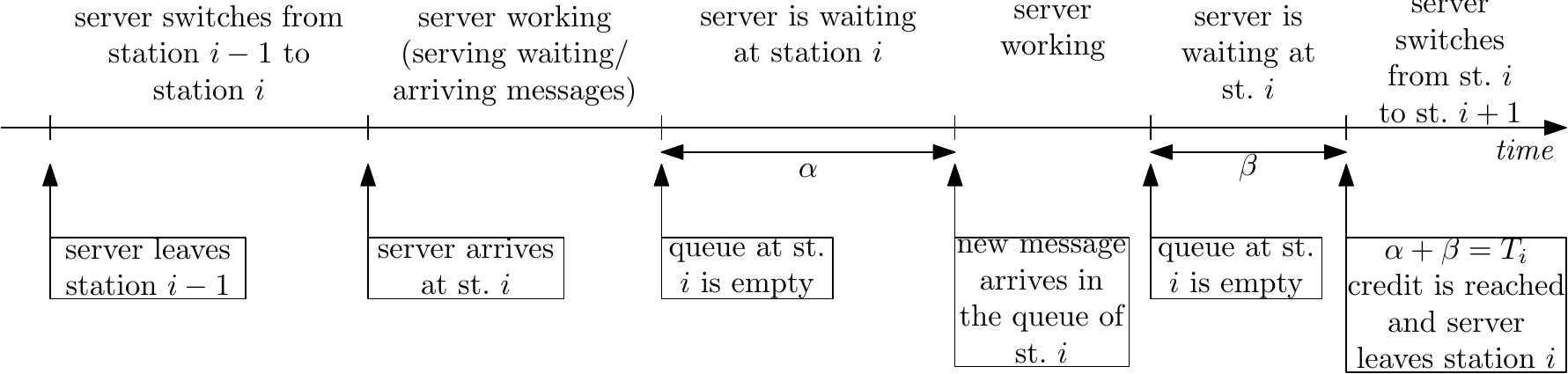} 
\caption{Operation of our polling model} \label{fig:operation}
\end{figure}

The message generation process, the lengths of the messages, and the switchover times are assumed to be independent -- both among each other and with respect to the other processes and stations. An illustration of the operation of the server is given in Figure~\ref{fig:operation}.

The goal of this paper is to derive an explicit formula for the mean average delay of a job for this model in steady state, that is, the expected time a message experiences from the point in time when it arrives in one of the queues until its service starts (i.e.\ excluding the processing time). The expected delay of a message generated at station $i$ is denoted by $\E D_i$, the mean average queueing delay is then defined by
$$
\bar{D}=\sum_{i=1}^N \frac{\rho_i}{\rho_0}\, \E D_i,
$$
where $\rho_i:=\lambda_i b_i$ is the traffic load offered to station $i$ and $\rho_0:=\sum_{i=1}^N \rho_i$ is the total load offered to the system. We stress that the delays of the different stations are weighted by the traffic intensity $\rho_i$, which implicitly includes weighting with the average message lengths, whereas the delays $\E D_i$ do not include weighting the delay of the individual messages with their lengths. This seems to be common in the literature; e.g.\ Takagi \cite{takagi} calls this quantity \emph{intensity weighted mean waiting time} (p.\ 92, \cite{takagi}).

The mean delay will be expressed explicitly in terms of the above parameters $\lambda_i$, $b_i$, $b_i^{(2)}$, $T_i$, $r_i$, and $r_i^{(2)}$, $i=1,\ldots, N$.

Furthermore, it will be convenient to use the following abbreviations.  We abbreviate by $r_0:=\sum_{i=1}^N r_i$ the sum of the mean switchover times and by $r_0^{(2)}:=\sum_{i=1}^N r_i^{(2)}+\sum_{i,j=1,i\neq j}^N r_i r_j$ the second moment of the sum of all switchover times. Finally, we let $T_0:=\sum_{i=1}^N T_i$ denote the total ``wait-and-see'' time per cycle.

\subsection{Main results}\label{sec:main}
In this section, we give our main results. Theorem~\ref{thm:main} gives a formula for the mean average delay in terms of the parameters of the system $\lambda_i$, $b_i$, $b_i^{(2)}$, $r_i$, and $r_i^{(2)}$, $i=1,\ldots, N$, as well as the times $T_i$, $i=1,\ldots, N$. This is simplified for the case of two stations, $N=2$, in Corollary~\ref{cor:twost}.

The formula for the delay allows to investigate the following question: Given the system parameters $\lambda_i$, $b_i$, $b_i^{(2)}$, $r_i$, and $r_i^{(2)}$, -- how does one have to adjust the parameters $T_i$, $i=1,\ldots, N$, such that the mean average delay is minimized. We will see that in many cases it is favourable -- in the sense of lower average queueing delay -- to choose \emph{positive} $T_i$. This is either due to (i) the random switchover times or (ii) the asymmetry of the system. This is described in detail for the case of two stations in Theorems~\ref{thm:optt1} and~\ref{thm:optt2}, where the effects (i) and (ii) are treated somehow in a decoupled way.

Finally, in Theorem~\ref{thm:lower} we consider all strategies that allow the server to wait at a station even though no work may be present. We give a lower bound for the delay for any such strategy.

Before we come to the main results, let us mention the stability condition for the system. Due to the exhaustive nature of our polling model it is clear that the system is stable if and only if
$$
\rho_0 < 1.
$$
We assume this condition from now on.

Furthermore, we recall the result for the so called exhaustive strategy from \cite{takagi}. In our model, this corresponds to the special case $T_1=T_2=\ldots=T_N=0$. In this case, one can find that
$$
\bar{D}=\frac{\sum_{i=1}^N\lambda_{i}b_{i}^{(2)}}{2(1-\rho_{0})}+\frac{r_{0}(\rho_{0}^{2}-\sum_{i=1}^N\rho_{i}^{2})}{2\rho_{0}(1-\rho_{0})}+\frac{r_{0}^{(2)}}{2r_{0}}.
$$
In terms of queueing delay, it was shown in \cite{luietal1992} that the exhaustive strategy provides the lowest delay in the class of all non-idle strategies, in particular, as compared to the gated and limited strategy. Therefore, it will serve as a benchmark for our strategy, which is a strategy allowing the server to be idle even though work may be present in the system (at other stations).

The main theorem for our polling model is as follows.
\begin{thm} \label{thm:main}
 The mean average delay of the polling model introduced above is given by:
\begin{eqnarray*}
\bar{D} & = & \frac{\sum_{i=1}^N\lambda_{i}b_{i}^{(2)}}{2(1-\rho_{0})}+\frac{(r_{0}+T_{0})(\rho_{0}^{2}-\sum_{i=1}^N\rho_{i}^{2})}{2\rho_{0}(1-\rho_{0})}+\frac{\frac{1}{2}\rho_{0}r_{0}^{(2)}+r_{0}\sum_{i=1}^N T_{i}(\rho_{0}-\rho_{i})}{\rho_{0}(r_{0}+T_{0})}\\
 &  & +\frac{1}{(r_{0}+T_{0})\rho_{0}} \left[ \sum_{i=1}^N T_{i}^2 \, \frac{(1-2\rho_{i})(\rho_{0}-\rho_{i})}{2(1-\rho_{i})}  + \sum_{1\leq i<j\leq N} T_i T_j (\rho_0-\rho_i-\rho_j) \right].
\end{eqnarray*}
\end{thm}

The proof of this theorem is given in Section~\ref{sec:general}. Certainly, one can ask which values of $T_1, \ldots, T_N$ lead to a minimal queueing delay. In other words, given the system parameters, we would like to know how we have to set $T_1, \ldots, T_N$ in order to minimize $\bar{D}$. Note that this is a non-trivial question, because the $T_i$ appear in numerator and denominator. In fact, this is a minimization problem in the variables $T_1, \ldots, T_N$, subject to the non-negativity restriction $T_i\geq 0$, for all $i=1,\ldots, N$, which can be carried out -- in principle -- explicitly. We discuss the respective minimizers below for $N=2$. Certainly, for large $N$, one would solve the problem numerically.

For two stations $N=2$, the main result reduces to the following simpler formula.
\begin{cor}\label{cor:twost}
The mean average delay of the polling model introduced above with $N=2$ is given by:
\begin{multline}
\bar{D} = \frac{\sum_{i=1}^2 \lambda_{i} b_{i}^{(2)}}{2(1-\rho_{0})} + \frac{1}{\rho_0(r_0+ T_0)}\left[\frac{r_0^{(2)}\rho_0}{2} +\frac{\rho_1 \rho_2 }{1-\rho_0}\,( r_0+T_0 )^2\right.\\ \left. +  \rho_2 T_1\left(r_0 + T_1\, \frac{1-2\rho_1}{2(1-\rho_1)}\right) + \rho_1 T_2\left(r_0 + T_2\, \frac{1-2\rho_2}{2(1-\rho_2)}\right)\right]. \label{eqn:formulaN2}
\end{multline}
\end{cor}

In particular, one can minimize (\ref{eqn:formulaN2}) w.r.t.\ $T_1$ and $T_2$ subject to the restrictions $T_1\geq 0$, $T_2\geq 0$ in order to obtain the minimal possible delay. Let us denote by $T_1^*$ and $T_2^*$ the minimizers. We say that there is \emph{no gain from waiting at station $i$} if $T_i^*=0$; if $T_i^*>0$ we say that it is \emph{worth waiting at station $i$}.

From the explicit expression above, one can observe the following consequences. First, we consider a partially symmetric polling model, by which we \emph{only} mean that both stations have the same intensities $\rho_1=\rho_2$, but not necessarily the same switchover time distribution, message length distribution, nor arrival rate (cf.\ \cite{pekoz1999}).

\begin{thm} \label{thm:optt1} Consider a polling model as introduced above with two stations. In the case of a symmetric polling model, $\rho_1=\rho_2$, the following holds.
 \begin{itemize}
  \item With deterministic switchover times, that is, $r_1^2=r_1^{(2)}$ and $r_2^2=r_2^{(2)}$, we get $T_1^*=T_2^*=0$. I.e.\ in this case there is no gain from waiting at either station.
  \item With non-deterministic switchover times, that is, $r_1^2<r_1^{(2)}$ or $r_2^2<r_2^{(2)}$, it is worth waiting (at both stations) if and only if
\begin{equation} \label{eqn:star}
2\rho_1 < 1 - \frac{r_0^2}{r_0^{(2)}+r_0^2 \,\frac{\rho_1}{1-2\rho_1}}.
\end{equation}
 In this case, the optimal waiting time $T_1^*=T_2^*>0$ can be calculated explicitly, see (\ref{eqn:optparforsymmetric}).  The minimal delay is then given by inserting (\ref{eqn:optparforsymmetric}) into (\ref{eqn:formulaN2}); this delay is strictly lower than the mean average delay induced by the exhaustive strategy.
 \end{itemize}
\end{thm}

We remark that the fraction on the right-hand side of (\ref{eqn:star}) equals 
\begin{equation}  \label{eqn:thefraction}
\frac{(\E[ R_1+R_2])^2}{ {\rm var}[R_1+R_2] + (\E[ R_1+R_2])^2 \, \frac{1-\rho_1}{1-2\rho_1} },
\end{equation}
where $R_i$ are independent switchover times for switching from station $i$ to station $i+1$, respectively.

Now we consider an asymmetric polling model, i.e.\ $\rho_1>\rho_2$.

\begin{thm} \label{thm:optt2} Consider a polling model as introduced above with two stations. In the case of an asymmetric polling model with deterministic switchover times, that is, assuming $\rho_1>\rho_2$ and $r_1^2=r_1^{(2)}$, $r_2^2=r_2^{(2)}$, the following holds.
 \begin{itemize}
 \item There is no gain from waiting at station 2, i.e.\ in all cases $T_2^*=0$.
 \item Further, it is worth waiting at station 1 if and only if
\begin{equation} \label{eqn:nor}
\rho_1-\rho_1^2 + \rho_2^ 2-\rho_2 - 2\rho_1 \rho_2 > 0.
\end{equation}
In this case, one can calculate the minimizer $T_1^*>0$ explicitly as in (\ref{eq:T1s}). The minimal delay is then given by inserting (\ref{eq:T1s}) and $T_2^*=0$ into (\ref{eqn:formulaN2}); this delay is strictly lower than the mean average delay induced by the exhaustive strategy.
 \end{itemize}
\end{thm}

\begin{rem}
A similar discussion is possible for the case of an asymmetric polling model with non-deterministic switchover times. There, both of the following effects will be combined. Namely, note that Theorem~\ref{thm:optt1} shows that large variances of the switchover times (increasing the variances of the switchover times in condition (\ref{eqn:star}), cf.\ (\ref{eqn:thefraction})) lead to the situation where it is worth waiting (at both stations). On the other hand, Theorem~\ref{thm:optt2} shows that a strong asymmetry (reducing $\rho_2$ in condition (\ref{eqn:nor})) makes it useful to wait at the station with significantly higher traffic intensity. These effects will both be present in the case of a not necessarily symmetric system with non-deterministic switchover times.
\end{rem}

\begin{rem}
We remark the following rather surprising fact: Note that (\ref{eqn:nor}) does not depend on the switchover times. So, the question whether it makes sense to wait at station $1$ only depends on the relation of the intensities $\rho_i$ and not on the length of the possible idle period due to the switching. Similarly, the expression in (\ref{eqn:star}) does not depend on the absolute lengths of the two switchover times (one can multiply both $R_i$ by the same constant without changing (\ref{eqn:star}), cf.\ (\ref{eqn:thefraction})) nor on the order of the switchover times (but only on the sum). However, even though the decision whether to wait or not does not depend on the absolute length of the switchover time, the resulting credit does (cf.\ (\ref{eqn:optparforsymmetric})).
\end{rem}

\begin{rem}
Similar discussions are possible for $N>2$ since $\bar{D}$ has the form
$$
\bar{D} = c + \frac{\vec{T} A \vec{T}^t + \vec{T} \vec{b}^t+a}{r_0+ T_0}.
$$
with some $N\times N$-matrix $A$, constants $c$, $a$, $\vec{b}$, and $\vec{T}:=(T_1,\ldots, T_N)$.
\end{rem}

\begin{rem}
 The trivial case $N=1$ is included in our results. It corresponds to a single queue where the server takes vacations. In our model, funnily enough, it takes a vacation (switchover time) after it has spent a total idle time of $T_1$ at the queue. The corresponding delay is
$$
\bar{D} = \frac{\lambda_{1} b_{1}^{(2)}}{2(1-\rho_{1})} + \frac{r_1^{(2)}}{2(r_1+ T_1)},
$$
which is obviously minimized for $T_1=\infty$, which also is easy to interpret, because then the server is never allowed to take a vacation; and the system is thus identical to an $M/G/1$ queue.
\end{rem}

Finally, we discuss a lower bound for the delay for strategies that allow waiting times of any type. We recall at this point that \cite{luietal1992} shows that the exhaustive strategy provides the lowest delay in the class of all non-idle strategies, that is, all strategies where the server is \emph{not} allowed to wait at a station if no work is present there.

In the following, we consider strategies that are not allowed to use future information of the system, that serve FCFS, and where the server is not idle if at its present station messages are waiting to be served. Further, we have to assume that with this strategy the system has a steady-state distribution.

The next theorem gives a lower bound for the delay for all of these strategies where the server \emph{is} allowed to wait at stations due to reasons that depend only on the current station in the current cycle (that is, since the server arrived at the present station). This restriction considers those strategies that look at the evolution of the traffic at the present station since the server arrived there. It does not allow strategies that take their decisions according to e.g.\ the queue status at different stations or the recent switchover times.

This provides a lower bound, in particular for the model treated so far, the strategy proposed in \cite{boxma2002}, as well as the strategy proposed in Section~\ref{sec:outlook} below.

\begin{thm} Consider a polling model where the stations are served in cyclic order. Then, for any strategy that allows the server to wait at a station even though no work is present there but the decision on whether and how long to wait only depends on the evolution of queue of the current station since the server arrived at the station, we have
\begin{eqnarray}
\bar{D} & \geq &  \frac{\sum_{i=1}^N \lambda_{i}b_{i}^{(2)}}{2(1-\rho_{0})} +  \frac{r_{0}(\rho_{0}^{2}-\sum_{i=1}^N\rho_{i}^{2})}{2\rho_0(1-\rho_{0})} \nonumber \\
 & & +\min_{f_1,\ldots, f_N\geq 0} \frac{1}{\rho_0(r_0+f_0)} \left[ \frac{\rho_0 r_0^{(2)}}{2}+ \sum_{j=1}^N \left(r_{0} f_{i}+\frac{f_{i}^{2}}{2}\right)(\rho_{0}-\rho_{i}) \right. \nonumber \\
 &  & \left.+\sum_{i=1}^Nf_{i}\left(\sum_{j=1}^{i-1}\alpha_j\left(\sum_{l=i+1}^N \rho_{l}+\sum_{l=1}^{j-1}\rho_{l}\right)+\sum_{j=i+1}^N \alpha_j\sum_{l=i+1}^{j-1}\rho_{l}\right) \right], 
\label{eq:Wschranke2}\end{eqnarray}
where $\alpha_j:=\rho_{j}\frac{r_{0}+f_{0}}{1-\rho_{0}}+f_{j}$ and $f_0=\sum_{i=1}^N f_i$. \label{thm:lower}
\end{thm}

The idea behind this theorem is that $f_i$ is the expected time the server spends waiting at station $i$ in a cycle. The time as such is random for a general strategy, of course. Since the $f_i$ are unknown in general, the minimum appears. In the case of the concrete model treated so far, we had $f_i=T_i$, because by definition the total time the server spends at station $i$ is deterministic and equals $T_i$.

The minimum in (\ref{eq:Wschranke2}) can be calculated explicitly in principle as well as numerically without any problem. A proof of Theorem~\ref{thm:lower} is given in Section~\ref{sec:lower}.

\subsection{Related work} \label{sec:relwork}
Basic references on polling models are \cite{takagi}, \cite{tag}, \cite{tag2}, \cite{luietal1992}.

References that refer to polling models where the server may be waiting at a station are apparently rare. The main references for us are Pek\"{o}z \cite{pekoz1999} and Boxma et al.\ \cite{boxma2002}. 

Pek\"{o}z introduced the strategy we use in this paper for the completely symmetric model (that is, all the arrival rates, service times, switchover times, and the $T_i$ are identical). In particular, his Theorem~2.2 is a special case of our Theorem~\ref{thm:main}. In the present paper, we consider the general polling model. Furthermore, a new observation is that also a sufficiently asymmetric system can make it useful to wait at a station, independently on whether the switchover times are random or not.

The second main reference is Boxma et al.\ \cite{boxma2002}, where a polling model with $N=2$ stations is analysed. In that work, the following situation is investigated: If the server encounters an \emph{empty queue} at station 1 \emph{once it arrives there}, a ``wait-and-see'' timer is activated in order to wait for the possible arrival of new messages. However, -- contrary to the present setup -- once the server has finished some work or the timer has run out, it will immediately switch to the next station. We compare the resulting delay obtained from this strategy to ours in Figure~2. We have found cases, where our strategy leads to lower delay than the strategy proposed by Boxma et al.\ and also cases where it performs worse. The latter is usually the case if the intensities $\rho_1, \rho_2$ are close to each other, whereas in the case that we deal with a highly asymmetric system, our strategy seems to be better. Also, for large switchover times our strategy seems to perform better than \cite{boxma2002}, since in this case the timer from \cite{boxma2002} is rarely activated. Unfortunately, it does not seem to be possible to compare the strategies directly due to the non-explicit nature of the delay formulas in \cite{boxma2002}.

 \addtocounter{figure}{1}
\begin{figure}
\begin{center}
\begin{tabular}[t]{ll}
\includegraphics[viewport=155 425 580 720, scale=0.5, clip=true]{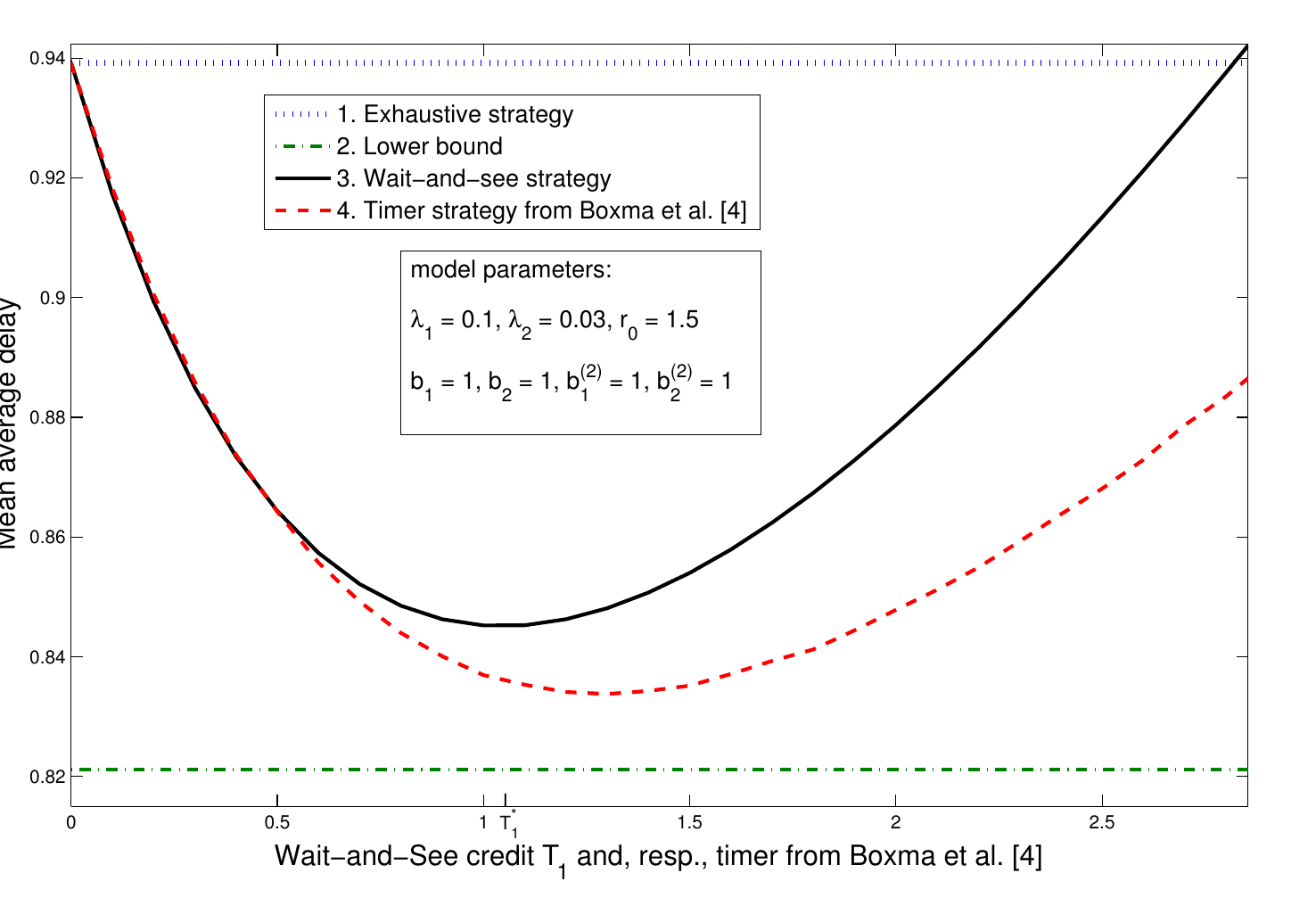}  & \includegraphics[viewport=155 425 580 720, scale=0.5, clip=true]{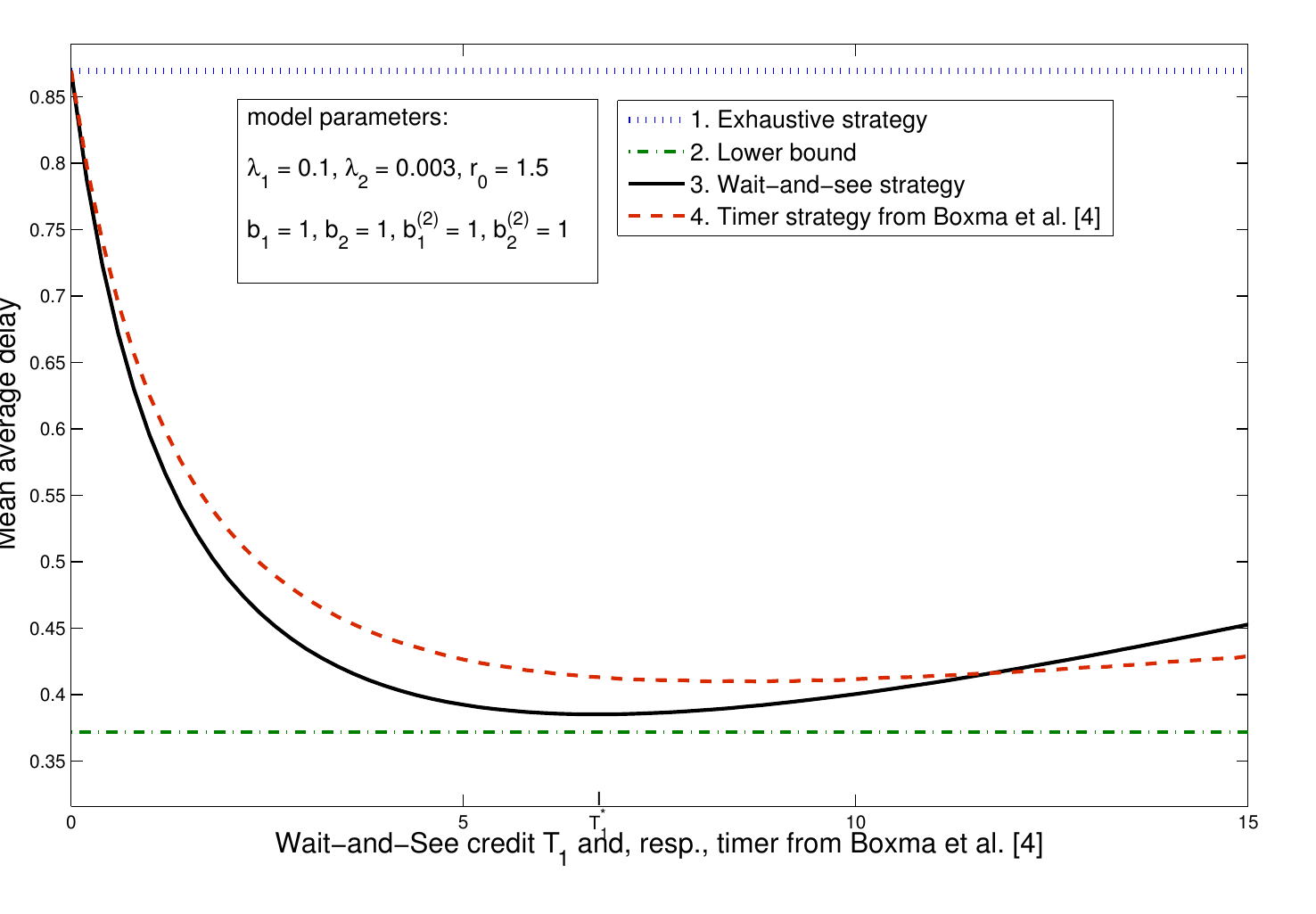} \\
\multirow{11}{*}{\includegraphics[viewport=155 425 580 720, scale=0.5, clip=true]{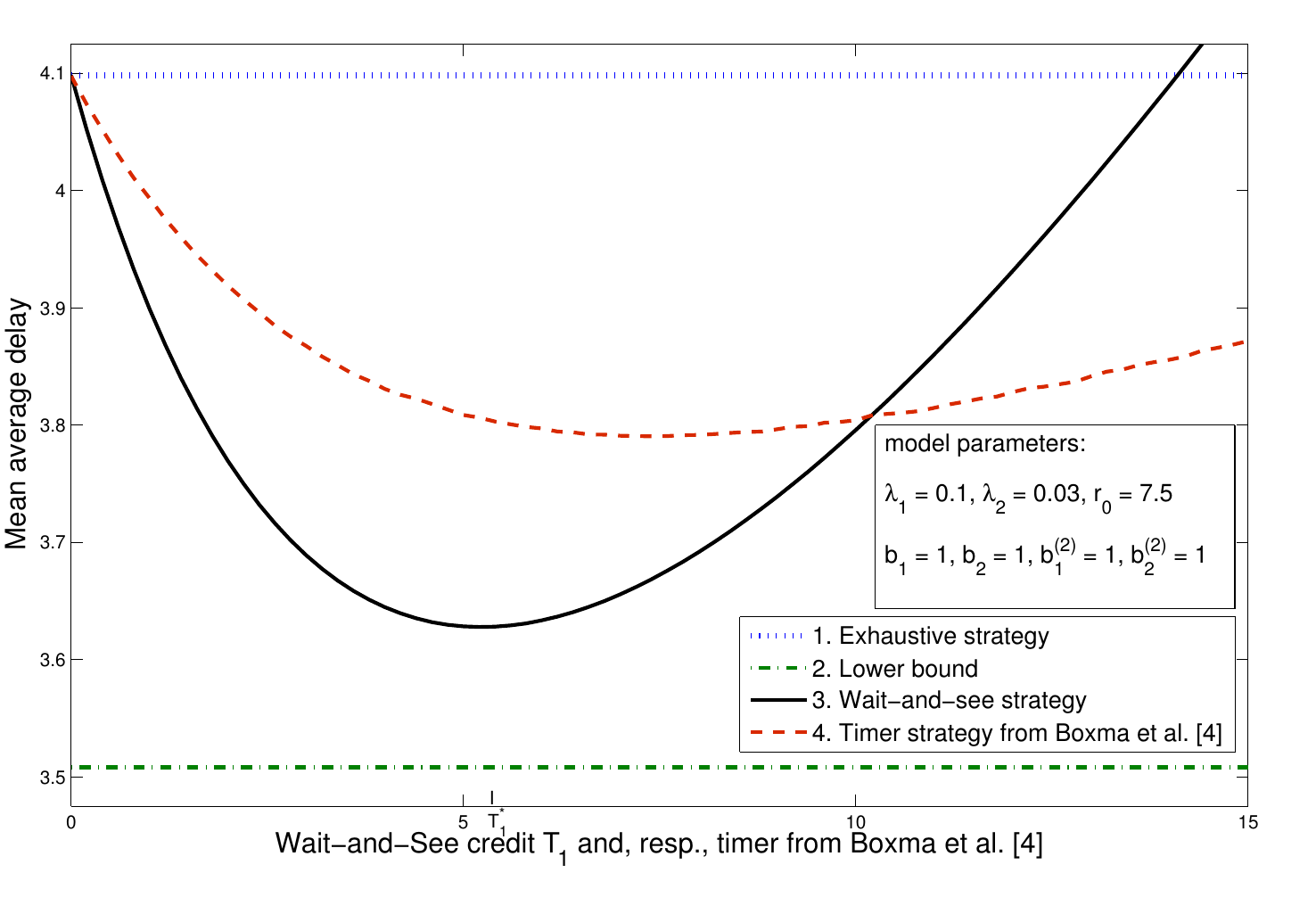} }
& Figure~\thefigure:\\
 & Comparison to Boxma et al.\ \cite{boxma2002},\\
  & delay vs.\ credit (resp.\ timer) at stat.\ 1.\\
  & Note that in the first plot the strategy \\
  & from \cite{boxma2002} provides lower delay; in the\\
  &  second plot only the arrival rate of the\\
  &   second station is changed and then our\\
  &   wait-and-see strategy has lower delay; \\
 &    in the third plot the same happens \\
 &  due to longer switchover times
\end{tabular}
\end{center}
\end{figure}

Further references on polling models where the server may be waiting consider only single-station systems with vacations (\cite{boxma2002b} and \cite{yechiali}).

The strategy employed in \cite{boxma2002} and in the present paper is somehow related to a so-called {\it forced idle time}. We refer e.g.\ to \cite{cooper1998,cooper1999} for some work on this. However, in the present setup, the server is not forced to be idle; whenever it is set to ``wait-and-see'', it rather resumes service as soon as new messages arrive. This is the reason we prefer the term ``wait-and-see'' rather than ``forced idle time''.

\section{The fundamental relations for the general polling model}\label{sec:general}
In this section, we derive the fundamental relations for the general polling model that allow us to obtain the formula for the mean average queueing delay. In particular, we give a proof of Theorem~\ref{thm:main}. We proceed in several steps.

\paragraph{Cycle time.} The cycle time is defined to be the time that the server takes from one arrival at station 1 to its next arrival at the same station. We obtain the average cycle time, $\E C$.

First, let us define more precisely the notion of the server being \emph{idle}, \emph{switching}, and \emph{waiting}. The server is \emph{waiting} when it is at some station waiting for messages to arrive. Note that, by the definition of our strategy, the total time the server spends waiting in each cycle equals $\sum_{i=1}^N T_i$. The server is said to be in the state of \emph{switching} from the time it leaves one station until it starts working at the next station. Finally, we say that the server is \emph{idle} if it is either waiting or switching.

Now, note that a cycle contains periods when the server works and periods when it is idle. In our polling model, the server is idle exactly for the time it waits and for the switchover time. Thus, the expected time the server is idle in a cycle equals
\begin{equation} \label{eqn:dhash}
(1-\rho_0) \E C  = \sum_{i=1}^N r_i + \sum_{i=1}^N T_i.
\end{equation}
This allows to deduce the expected cycle time in our polling model in steady state:
\begin{equation} \label{eqn:hash}
\E C = \frac{r_0 +  T_0}{1-\rho_0}.
\end{equation}

\paragraph{The decomposition principle.} We will use a decomposition principle to deduce our results. One can find a proof of this principle in other contexts e.g.\ in \cite{groenendijk,boxma2002}. We omit the proof for our system since it is completely analogous.

In order to formulate the decomposition principle, we need the notion of the \emph{workload} $V$, which we define to be the sum of all message lengths that are present in the system (including the remaining length of the currently processed message) at a random point in time in steady state.

Then the decomposition principle says that $V$ has the same distribution as
$$
V \deq V^{\rm{M/G/1}}+V^{\rm{idle}},
$$
where $V^{\rm{M/G/1}}$ is the workload in the same polling model with no switching or waiting times, that is, an M/G/1 queue. On the other hand, $V^{\rm{idle}}$ is the workload at a random point in time given that the server is idle at that point, and $V^{\rm{M/G/1}}$ and $V^{\rm{idle}}$ are independent.

Let
$$
q:=\pr{\text{server switching}~|~\text{server idle}} = \frac{\pr{\text{server switching}}}{\pr{\text{server switching}}+\pr{\text{server waiting}}}.
$$
Therefore,
\begin{equation} \label{eqn:decomp}
\E V = \E V^{\rm{M/G/1}}+q \E V^{\rm{switching}}+ (1-q) \E V^{\rm{waiting}},
\end{equation}
where $V^{\rm{switching}}$ and $V^{\rm{waiting}}$ are the workloads at a random point in time given that the server is switching and waiting, respectively, at that point.

\paragraph{Expected workload.}
We now calculate $\E V$ in two different ways. On the one hand, note that
\begin{equation} \label{eqn:workload1}
\E V =  \sum_{i=1}^N b_{i}\E[\#\, \text{messages in queue at station $i$}]+\sum_{i=1}^N\rho_{i}\frac{b_{i}^{(2)}}{2b_{i}}.
\end{equation}
Indeed, the first term accounts for the fact that there are messages that are not yet in service and that are waiting at the different stations. The second term corresponds to the fact that with probability $\rho_i$ we are looking at station $i$ and a message is being processed there. The workload of that message is exactly its expected residual lifetime, that is, $b_{i}^{(2)}/(2b_{i})$.

By Little's law, the last equation becomes
\begin{equation} \label{eqn:littleslaw}
\E V =  \sum_{i=1}^N b_{i} \lambda_i \E D_i +\sum_{i=1}^N\rho_{i}\frac{b_{i}^{(2)}}{2b_{i}} = \sum_{i=1}^N \rho_i \E D_i +\sum_{i=1}^N\rho_{i}\frac{b_{i}^{(2)}}{2b_{i}} = \rho_0 \bar{D} + \sum_{i=1}^N\rho_{i}\frac{b_{i}^{(2)}}{2b_{i}}.
\end{equation}
This equation shows that, in order to obtain the mean average queueing delay of our system, $\bar{D}$, we have to calculate the expected workload $\E V$.

On the other hand, we use the decomposition principle (\ref{eqn:decomp}). Clearly, $\E V^{\rm{M/G/1}}$ is known:
$$
\E V^{\rm{M/G/1}} = \frac{\sum_{i=1}^N \lambda_{i}b_{i}^{(2)}}{2(1-\rho_{0})},
$$
see e.g.\ \cite{kleinrock}, p.\ 201. Therefore, we obtain with (\ref{eqn:decomp}) and (\ref{eqn:littleslaw}) that
\begin{equation} \label{eqn:wlmain}
\rho_0 \bar{D} =  \frac{\sum_{i=1}^N \lambda_{i}b_{i}^{(2)}}{2(1-\rho_{0})}+q \E V^{\rm{switching}}+ (1-q) \E V^{\rm{waiting}} - \sum_{i=1}^N\rho_{i}\frac{b_{i}^{(2)}}{2b_{i}}.
\end{equation}

We are left with calculating the expected workload given that we find the system in the state of switching and, respectively, the expected workload given that we find the system in the state of waiting, as well as some relation between time periods of switching and waiting. As concerns the latter, in fact, it is sufficient that we clearly know that
\begin{equation} \label{eqn:relationbetweenq1}
 \frac{q}{\pr{\text{server switching}}}= \frac{1}{\pr{\text{server switching}}+\pr{\text{server waiting}}}=\frac{1}{1-\rho_0}
\end{equation}
and
\begin{equation} \label{eqn:relationbetweenq2}
 \frac{1-q}{\pr{\text{server waiting}}}= \frac{1}{\pr{\text{server switching}}+\pr{\text{server waiting}}}=\frac{1}{1-\rho_0}.
\end{equation}

\paragraph{Workload present while switching.} Observe that
\begin{equation} \label{eqn:wlmain2}
\E V^{\rm{switching}} = \frac{1}{\pr{\text{server switching}}}\, \sum_{i=1}^N p_i \E V_i^{\rm{switching}},
\end{equation}
where $p_i$ is the probability of encountering the server in the state of switching from station $i$ to station $i+1$ when entering the system at a random point in time and $\E V_i^{\rm{switching}}$ is the expected workload at such a point in time.

Clearly, $p_i=r_i/\E C$, since this is the fraction of time in a cycle that the server spends switching from station $i$ to station $i+1$.

Now, we have to find $\E V_i^{\rm{switching}}$, the expected total amount of work that is present given that we look at the system at a point when the server is switching from station $i$ to station $i+1$. Let us assume we are at such a point in time; then there are different times when the currently present workload was generated. We distingish these times and determine the respective workload:

\begin{itemize}
 \item At all stations workload was generated during the current switching period. It is given by the expected residual lifetime (in fact, the expected backwards recurrence time) of the current switching period: $\frac{r_i^{(2)}}{2 r_i}\cdot\sum_{j=1}^N \rho_j $.
 \item At all stations, except for the $i$-th station, workload was generated while the server was at station $i$ (working and waiting).
The time spent working has mean $\rho_i \E C$, and the time spent waiting equal $T_i$, so that the generated workload becomes $(\rho_i \E C+T_i)\cdot\sum_{j=1, j\neq i}^N \rho_j$.
 \item Similarly, while the server was at station $k$ (working and waiting), $k\neq i$, at all stations (except for those that later in the cycle became emptied) workload was generated and is still present. The time spent working has mean $\rho_k \E C$, the time spent waiting equals $T_k$, so that the generated workload becomes
\begin{equation} \label{eqn:wrkldworking}
(\rho_k \E C +T_k) \cdot \sum_{j\in\{i+1,\ldots, k-1\}} \rho_j,
\end{equation}
where $\{i+1,\ldots, k-1\}$ is defined as $\{i+1,\ldots, N\}\cup \{1,\ldots, k-1\}$ if $i+1>k-1$.
 \item During the switchover time from station $k$ to station $k+1$ ($k\neq i$), which takes on average $r_k$, workload was generated at all stations (except for those that later in the cycle became emptied): 
\begin{equation} \label{eqn:wrkldswitching}
r_k \cdot \sum_{j\in\{i+1,\ldots, k\}} \rho_j.
\end{equation}
\end{itemize}
Summing up all this workload, we get
\begin{eqnarray}
\E V_i^{\rm{switching}} & = & \sum_{j<i} r_{j}\left(\sum_{l=i+1}^{N}\rho_{l}+\sum_{l=1}^{j}\rho_{l}\right)+\sum_{j>i} r_{j}\sum_{l=i+1}^j\rho_{l}\nonumber \\
 &  & +\sum_{j<i} \rho_{j}\mathbb{E}C\left(\sum_{l=i+1}^N\rho_{l}+\sum_{l=1}^{j-1}\rho_{l}\right)+\sum_{j>i}\rho_{j}\mathbb{E}C\sum_{l=i+1}^{j-1}\rho_{l}\nonumber \\
 &  & +\sum_{j<i} T_{j}\left(\sum_{l=i+1}^N\rho_{l}+\sum_{l=1}^{j-1}\rho_{l}\right)+\sum_{j>i} T_{j}\sum_{l=i+1}^{j-1}\rho_{l}\nonumber \\
 &  & +\rho_i\E C (\rho_0-\rho_i)+(\rho_{0}-\rho_{i})\, T_{i}+\rho_{0}\frac{r_i^{(2)}}{2r_{i}}.\label{eq:EVi_schaltet}\end{eqnarray}

\paragraph{Workload present while waiting.} Analogously to the workload while switching, we observe that
\begin{equation} \label{eqn:wlmain3}
\E V^{\rm{waiting}} = \frac{1}{\pr{\text{server waiting}}}\, \sum_{i=1}^N q_i \E V_i^{\rm{waiting}},
\end{equation}
where $q_i$ is the probability of finding the server waiting (in a wait-and-see state) for messages at station $i$ and $\E V_i^{\rm{waiting}}$ is the workload one would find at such a point in time.

Clearly, $q_i=T_i/\E C$, since this is the fraction of time in a cycle that the server spends waiting at station $i$, by the definition of our polling model.

Similarly to the workload while switching, we obtain the expected workload generated while station $i$ is in the state of waiting. Let us assume we are at such a point in time, then there are different times when the currently present workload was generated. We distinguish these times and determine the respective workload:
\begin{itemize}
 \item When the server started working at station $i$, there was work waiting there. We denote the length of the ``busy period'' generated by this waiting traffic by $Z_i$. The workload generated at the other stations during this busy period is
\begin{equation} \label{eqn:terma}
\E Z_i \cdot (\rho_0 -\rho_i).
\end{equation}
In order to determine $\E Z_i$ note that the average time the server spends working at station $i$ is, on the one hand, $\rho_i \E C$. On the other hand, the time the server spends working at station $i$ consists of the length of the first busy period, $\E Z_i$, and all other busy periods generated, which is the number of busy periods in an M/G/1 queue with total idle time $T_i$. However, the expected number of busy periods in an M/G/1 queue with total idle time $T_i$ is $\lambda_i T_i$ (just disregard the time of the busy periods, then because of the memoryless property, the number of busy periods is Poisson with intensity $\lambda_i T_i$).

Thus,
$$
\rho_i \E C = \E Z_i + \lambda_i T_i\, \frac{b_i}{1-\rho_i},
$$
since $\frac{b_i}{1-\rho_i}$ is the average length of one busy period at station $i$. This equation allows to calculate $\E Z_i$.
 \item Workload was generated at all other stations except for the $i$-th, during the busy periods that have already taken place at station $i$, not considering the first busy period when the server started to work at station $i$. Per busy period, a workload at the other stations of in total
$$
\frac{b_i}{1-\rho_i} \cdot (\rho_0-\rho_i)
$$
was generated, since $\frac{b_i}{1-\rho_i}$ is the average length of one busy period at station $i$. In order to obtain the number of busy periods that have already taken place, note that these are on average $\lambda_i T_i/2$, because the waiting time is deterministic. Thus, we obtain
\begin{equation} \label{eqn:termb}
\frac{\lambda_i T_i}{2}\, \frac{b_i}{1-\rho_i} \cdot (\rho_0-\rho_i)
\end{equation}
for the total workload that was generated at all other stations during the busy periods (except for the very first one) at station $i$. 
 \item The total workload that was generated at all other stations during the waiting time spent so far at station $i$ is with the same reasoning
\begin{equation} \label{eqn:termc}
\frac{T_i}{2} \cdot (\rho_0-\rho_i).
\end{equation}
 \item The term in (\ref{eqn:wrkldworking}) has to be considered in the same way.
 \item The term in (\ref{eqn:wrkldswitching}) has to be considered in the same way.
\end{itemize}
Summing up all this workload, we get
\begin{eqnarray}
\E V_i^{\rm{waiting}} & = & \sum_{j<i}r_{j}\left(\sum_{l=i+1}^N\rho_{l}+\sum_{l=1}^j \rho_{l}\right)+\sum_{j>i}r_{j}\sum_{l=i+1}^{j}\rho_{l}\nonumber \\
 &  & +\sum_{j<i}\rho_{j}\mathbb{E}C\left(\sum_{l=i+1}^N\rho_{l}+\sum_{l=1}^{j-1}\rho_{l}\right)+\sum_{j>i}\rho_{j}\mathbb{E}C \sum_{l=i+1}^{j-1}\rho_{l}\nonumber \\
 &  & +\sum_{j<i}T_{j}\left(\sum_{l=i+1}^N\rho_{l}+\sum_{l=1}^{j-1}\rho_{l}\right)+\sum_{j>i}T_{j}\sum_{l=i+1}^{j-1}\rho_{l}\nonumber \\
 &  & +(\rho_{0}-\rho_{i})\left(\rho_{i}\mathbb{E}C+\frac{T_{i}}{2}(1-\frac{\rho_{i}}{1-\rho_{i}})\right).\label{eq:EVi_wartet}\end{eqnarray}

\begin{proof}[ of Theorem~\ref{thm:main}] In order to see the formula in Theorem~\ref{thm:main}, one just has to combine (\ref{eqn:wlmain}), (\ref{eqn:hash}), and (\ref{eqn:relationbetweenq1}), (\ref{eqn:wlmain2}), (\ref{eq:EVi_schaltet}), and (\ref{eqn:relationbetweenq2}), (\ref{eqn:wlmain3}), (\ref{eq:EVi_wartet}).\end{proof}

\section{The case of two stations}\label{sec:two}
In this section, we prove Theorems~\ref{thm:optt1} and~\ref{thm:optt2}. First we prove that the optimal parameters in the two-station case satisfy a linear relation. Then we prove Theorem~\ref{thm:optt1} (symmetric case, deterministic and non-deterministic) and Theorem~\ref{thm:optt2} (asymmetric and deterministic), respectively.

For simplicity, we introduce the following abbreviations:
\begin{eqnarray*}
c_{1} & := & \frac{\sum_{i=1}^2\lambda_{i}b_{i}^{(2)}}{2(1-\rho_{0})}\\
c_{2} & := & \frac{\rho_{1}\rho_{2}r_{0}^{2}}{1-\rho_{0}}+\frac{\rho_{0} r_0^{(2)}}{2}\\
c_{3} & := & r_{0}\rho_{2}+\frac{2 \rho_{2}\rho_{1}r_{0}}{1-\rho_{0}}\\
c_{4} & := & r_{0}\rho_{1}+\frac{2 \rho_{1}\rho_{2}r_{0}}{1-\rho_{0}}\\
c_{5} & := & \frac{2\rho_{2}\rho_{1}}{1-\rho_{0}}\\
c_{6} & := & \frac{c_5}{2}+\frac{\rho_{2}}{2}\left(1-\frac{\rho_{1}}{1-\rho_{1}}\right)\\
c_{7} & := & \frac{c_5}{2}+\frac{\rho_{1}}{2}\left(1-\frac{\rho_{2}}{1-\rho_{2}}\right).\end{eqnarray*}

An easy calculation shows that these are non-negative constants. With these abbreviations, formula (\ref{eqn:formulaN2}) becomes:
\begin{equation}\label{eqn:formulaN2withcs}
 \bar{D}= c_{1}+\frac{c_{2}+c_{3}T_{1}+c_{4}T_{2}+c_{5}T_{1}T_{2}+c_{6}T_{1}^{2}+c_{7}T_{2}^{2}}{\rho_{0}(r_{0}+T_{1}+T_{2})}.
\end{equation}

\begin{lem} On the set $r_0+T_0>0$, the minimizers of the quantity in (\ref{eqn:formulaN2withcs}) satisfy the following linear relation:
\begin{equation} 
(c_{5}-2c_{6})T_{1}^{*}=c_{3}-c_{4}+(c_{5}-2c_{7})T_{2}^{*}.\label{eq:LinZusam}
\end{equation}
In particular, in the symmetric polling model ($\rho_1=\rho_2$), we must have
\begin{equation}
T_{1}^{*}=T_{2}^{*}.\label{eq:T1_T2}\end{equation}\label{lem:lem}
\end{lem}

We remark that the above minimizers can be negative. Recall that we are interested in the optimal parameter, which are the minimizers of (\ref{eqn:formulaN2withcs}) subject to the restriction $T_1^*,T_2^*\geq 0$. This is why we distinguish in the following between the minimizers of (\ref{eqn:formulaN2withcs}) and the optimal parameters.

\begin{proof}[ of Lemma~\ref{lem:lem}]
Clearly, (\ref{eqn:formulaN2withcs}) shows that $\bar{D}$ can be written as follows:
\[\bar{D}=\bar{D}(T_1,T_2)=c_{1}+\frac{f(T_{1},T_{2})}{\rho_{0}(r_{0}+T_{1}+T_{2})},\]
with some function $f$. If $\bar{D}$ has a minimum at $T_1^*$ and $T_2^*$ (with $r_0+T_1^*+T_2^*>0$) it must satisfy
$$
\frac{\partial \bar{D}}{\partial T_1}(T_1^*,T_2^*)=0\qquad\text{and}\qquad\frac{\partial \bar{D}}{\partial T_2}(T_1^*,T_2^*)=0.
$$
Due to the fact that the denominator is a linear function in $T_1+T_2$, an easy calculation shows that we must actually have
$$
\frac{\partial f}{\partial T_1}(T_1^*,T_2^*)=\frac{\partial f}{\partial T_2}(T_1^*,T_2^*).
$$
This is
\[
c_{3}+c_{5}T_{2}^{*}+2c_{6}T_{1}^{*}=c_{4}+c_{5}T_{1}^{*}+2c_{7}T_{2}^{*},\]
exactly as asserted in (\ref{eq:LinZusam}).

In the symmetric case we have $\rho_1=\rho_2<1/2$ and $c_3=c_4$ and $c_7=c_6$ which implies that (\ref{eq:LinZusam}) becomes (\ref{eq:T1_T2}).
\end{proof}

\paragraph{Symmetric polling model.} We now consider a symmetric polling model, i.e.\ $\rho_1=\rho_2=:\rho$.

\begin{proof}[ of Theorem~\ref{thm:optt1}]
Assume that $T_1^*>0$ and $T_2^*>0$ are the optimal parameters. Then we know from (\ref{eq:T1_T2}) that $T_{1}^{*}=T_{2}^{*}=:T$. Therefore, we obtain:
\begin{eqnarray*}
\bar{D} 
 & = & c_{1}+\frac{(r_{0}^{(2)}+\rho\frac{r_{0}^{2}}{1-2\rho})+2(r_{0}+\frac{2\rho r_{0}}{1-2\rho})T+(\frac{4\rho}{1-2\rho}+1-\frac{\rho}{1-\rho})T^{2}}{2(r_{0}+2T)}.\end{eqnarray*}
The minimum of this expression is attained at
\begin{equation} \label{eqn:optparforsymmetric}
T^{*}=-\frac{1}{2}\,r_{0}+\frac{1}{2}\,\sqrt{4r_{0}^{2}\rho-3r_{0}^{2}+\left(r_{0}^{(2)}+r_{0}^{2}\frac{\rho}{1-2\rho}\right)\left(4-12\rho+8\rho^{2}\right)}.
\end{equation}
Let $a:=r_{0}^{(2)}+r_{0}^{2}\frac{\rho}{1-2\rho}$. The condition for $T^*$ to be well-defined and positive is:
$$
4r_{0}^{2}\rho-3r_{0}^{2}+a \cdot 4 (1-\rho)(1-2\rho) > r_0^2.
$$
This is true if and only if
$$
a(1-2\rho)>r_0^2,
$$
which is easily seen to be equivalent to what we stated in (\ref{eqn:star}).

In the deterministic case, $a=r_{0}^{2}+r_{0}^{2}\frac{\rho}{1-2\rho}$,
the condition becomes
$$
2\rho<1-\frac{r_{0}^{2}}{r_{0}^{2}+r_{0}^{2}\frac{\rho}{1-2\rho}}=1-\frac{1}{1+\frac{\rho}{1-2\rho}},
$$
which can easily be seen to lead to the contradiction $\rho<0$.
\end{proof}

\paragraph{Unsymmetric polling model with deterministic switchover times.} We now consider an a asymmetric polling model, i.e.\ $\rho_1>\rho_2$ with deterministic switchover times, i.e.\ $r_1^{(2)}=r_1^2$ and $r_1^{(2)}=r_1^2$.

\begin{proof}[ of Theorem~\ref{thm:optt2}, first part]
Recall that we would like to show that there is no gain from waiting at the station with less traffic, that is, station 2 in our case. We distinguish two cases: $\rho_1>1/2$ and $\rho_1<1/2$.

\emph{First case:} $\rho_{1}>\rho_{2}$ and $\rho_{1}>1/2$.

Note that trivially $\rho_{2}<1/2$. Recall that the linear relation (\ref{eq:LinZusam}) holds for the minimizers of (\ref{eqn:formulaN2withcs}) (which is the same as (\ref{eqn:formulaN2})). Since $\rho_{1}>1/2$ and $\rho_{2}<1/2$ we get $\frac{c_{3}-c_{4}}{c_{5}-2c_{6}}<0$ and $\frac{c_{5}-2c_{7}}{c_{5}-2c_{6}}<0$. Therefore, due to (\ref{eq:LinZusam}) one of the minimizers $T_1^*$ or $T_2^*$ must be negative. Therefore, the minimizers subject to the restriction $T_1^*, T_2^*\geq 0$ must satisfy either $T_{2}^{*}=0$ or $T_{1}^{*}=0$. However, the second case can be excluded easily: If we set $T_{1}=0$ in (\ref{eqn:formulaN2withcs}) and optimize in $T_{2}$ we would get:
$$
T_{2}^*=-r_{0}+\sqrt{r_{0}^{2}+\frac{c_2-c_{4}r_{0}}{c_7}}.
$$
This can be seen to be negative, because
$$
c_{2}-c_{4}r_{0}<0
$$
follows from
$$
-r_{0}^{2}\rho_{1}-\frac{\rho_{1}\rho_{2}r_{0}^{2}}{1-\rho_{0}}+\frac{1}{2}\rho_{0}r_{0}^{(2)}=r_{0}^{2}\frac{1}{2}(\rho_{2}-\rho_{1})-\frac{\rho_{1}\rho_{2}r_{0}^{2}}{1-\rho_{0}}<0,
$$
which holds since $\rho_2<\rho_1$.

Therefore, the case $T_{1}^{*}=0$, $T_{2}^{*}>0$ can be excluded; and we must have $T_{1}^{*}\geq 0$ and $T_{2}^{*}=0$ for the minimizers of (\ref{eqn:formulaN2withcs}) subject to $T_1^*, T_2^*\geq 0$.

\emph{Second case:} $\rho_{1}>\rho_{2}$ and $\rho_{1}<1/2$.

First let us rewrite the delay formula (\ref{eqn:formulaN2}). We exclude the trivial case $r_0=0$ and set $S_{1}:=T_{1}/r_{0}$ and $S_{2}:=T_{2}/r_{0}$. Then (\ref{eqn:formulaN2}) becomes:
\begin{equation} \label{eqn:delayintermsofS}
 \bar{D}=c_{1}+\frac{r_{0}}{\rho_{0}}\left(\frac{\rho_{1}\rho_{2}}{1-\rho_{0}}(1+S_{1}+S_{2})+\frac{\frac{1}{2}\rho_{0}+S_{1}\rho_{2}+S_{2}\rho_{1}+S_{1}^{2}\frac{\rho_{2}}{2}\frac{1-2\rho_{1}}{1-\rho_{1}}+S_{2}^{2}\frac{\rho_{1}}{2}\frac{1-2\rho_{2}}{1-\rho_{2}}}{1+S_{1}+S_{2}}\right).
\end{equation}

With the notation $S_{1}$ and $S_{2}$, the linear relation (\ref{eq:LinZusam}) becomes
$$
S_{1}=\frac{c_{3}-c_{4}}{c_{5}-2c_{6}}\frac{1}{r_{0}}+\frac{c_{5}-2c_{7}}{c_{5}-2c_{6}}S_{2}=\frac{\rho_{1}-\rho_{2}}{\rho_{2}\frac{1-2\rho_{1}}{1-\rho_{1}}}+\frac{\rho_{1}\frac{1-2\rho_{2}}{1-\rho_{2}}}{\rho_{2}\frac{1-2\rho_{1}}{1-\rho_{1}}}S_{2}.
$$
Setting $c:=\frac{\rho_{1}-\rho_{2}}{\rho_{2}\frac{1-2\rho_{1}}{1-\rho_{1}}}$
and $b:=\frac{\rho_{1}\frac{1-2\rho_{2}}{1-\rho_{2}}}{\rho_{2}\frac{1-2\rho_{1}}{1-\rho_{1}}}$, this is
\begin{equation}S_{1}=c+bS_{2}.\label{eq:S1}\end{equation}

Consider $\bar{D}=\bar{D}(S_{1},S_{2})$ (given in (\ref{eqn:delayintermsofS})) as a function of $S_1$ and $S_2$. It suffices to consider the function $\bar{D}(S_{1},S_{2})$ only at those points that satisfy the linear relation (\ref{eq:S1}), that is, $\bar{D}(c+b S_{2},S_{2})$, $S_{2}\in[0,\infty)$. We are finished if we can show that the derivative of this function w.r.t.\ $S_2$ at $0$ is non-negative, since then the optimum must be attained for negative $S_2$, and thus negative $T_2$, which is impossible. Let $g(S_{2}):=\bar{D}(c+b S_{2},S_{2})$.

The function $g$ can be written as follows:
\begin{eqnarray*}
g(S_{2}) & = & c_{1}+\frac{r_{0}}{\rho_{0}}\left(\frac{\rho_{1}\rho_{2}}{1-\rho_{0}}(1+c+b S_{2}+S_{2})\right.\\
 &  &\left. +\frac{\frac{1}{2}\rho_{0}+(c+bS_{2})\rho_{2}+S_{2}\rho_{1}+(c+bS_{2})^{2}\frac{\rho_{2}}{2}\frac{1-2\rho_{1}}{1-\rho_{1}}+S_{2}^{2}\frac{\rho_{1}}{2}\frac{1-2\rho_{2}}{1-\rho_{2}}}{1+c+bS_{2}+S_{2}}\right).\end{eqnarray*}
Clearly,
\begin{eqnarray*}
g'(0) & = & \frac{r_{0}}{\rho_{0}}\left(\frac{\rho_{1}\rho_{2}}{1-\rho_{0}}(1+b)\right.\\
 &  & \left. +\frac{(1+c)(\rho_{2}b+\rho_{1}+c\rho_{2}\frac{1-2\rho_{1}}{1-\rho_{1}})-(1+b)(\frac{1}{2}\rho_{0}+\rho_{2}c+c^{2}\frac{\rho_{2}}{2}\frac{1-2\rho_{1}}{1-\rho_{1}})}{(1+c)^2}\right).
\end{eqnarray*}
We would like to show that $g'(0)>0$, which is true if and only if
$$
\frac{\rho_{1}\rho_{2}(1+b)(1+c)^{2}}{1-\rho_{0}}+(\rho_{2}b+\rho_{1}+c\rho_{2}\frac{1-2\rho_{1}}{1-\rho_{1}}b)(1+c)-(\frac{\rho_{0}}{2}+\rho_{2}c+\frac{c^{2}\rho_{2}}{2}\frac{1-2\rho_{1}}{1-\rho_{1}})(1+b)>0
$$
After some calculations, it can be seen that this is equivalent to
$$
-2\rho_{1}\rho_{2}-\rho_{1}+\rho_{1}^{2}+\rho_{2}-\rho_{2}^{2}+2\rho_{1}\rho_{2}^{2}+2\rho_{1}^{2}\rho_{2}<0
$$
which is easily seen to be always satisfied in the case $\rho_1,\rho_2<1/2$.
\end{proof}

We have seen that in the symmetric polling model with deterministic switchover times there is no gain from waiting at the station with less traffic (station 2). Now, we determine when it is useful to wait at the station with more traffic (station 1), and what is the optimal waiting time $T_1^*$ in this case. It turns out that the condition is (\ref{eqn:nor}); and the optimal parameter is given by:
\begin{equation}
T_{1}^{*}=-r_{0}+\sqrt{r_{0}^{2}+\frac{c_2-c_{3}r_{0}}{c_6}}; \label{eq:T1s}\end{equation}
and the corresponding delay is then obtained by plugging in (\ref{eq:T1s}) and $T_2^*=0$ into (\ref{eqn:formulaN2withcs}).

\begin{proof}[ of Theorem~\ref{thm:optt2}, second part]
We get the optimal parameter if we set $T_2=0$ in (\ref{eqn:formulaN2withcs}) and differentiate w.r.t.\ $T_1$. Then the minimizer is seen to be given by (\ref{eq:T1s}). Condition (\ref{eqn:nor}) corresponds to $T_{1}^{*}>0$: in order for (\ref{eq:T1s}) to be positive, we must have
$$
c_{2}-c_{3}r_{0}>0,
$$
which translates into
$$
\frac{\rho_{0}}{2}>\rho_{2}+\frac{\rho_{2}\rho_{1}}{1-\rho_{0}},
$$
and thus (\ref{eqn:nor}) appears.
\end{proof}

\section{Lower bound}\label{sec:lower}
The goal of this section is to give a proof of Theorem~\ref{thm:lower}. For this purpose, let us define the following random variables. We denote by $F_i$ the time in steady state that the server spends at station $i$ waiting in a cycle, that is, being idle because there is no work at that station. Furthermore, let $f_i:=\E F_i$, $f_i^{(2)}:=\E F_i^2$, and $f_0:=\sum_{i=1}^N f_i$.

For a general strategy basically nothing can be said about the distribution of the $F_i$ even about their means $f_i$. The idea behind the proof of Theorem~\ref{thm:lower} is to estimate the mean average delay by an expression in terms of $f_i$ and $f_i^{(2)}$, to estimate by Jensen's inequality,
\begin{equation} \label{eqn:jensen}
 f_i^{(2)}\geq f_i^2,
\end{equation}
and thus to obtain an expression that only depends on the $f_i$. Then, minimizing over $f_i\geq 0$, we obtain the lower bound (\ref{eq:Wschranke2}). The details are as follows.

\paragraph*{Decomposition.} Let $V$ be the workload as defined above. Analogously to the decomposition principle in (\ref{eqn:decomp}) one can show that
$$
\E V = \E V^{\rm{M/G/1}}+\E V^{{\rm idle}}.
$$
We note that in the derivation of (\ref{eqn:workload1}) and (\ref{eqn:littleslaw}), the strategy was not used at all. So, one obtains a lower bound for $\bar{D}$ from a lower bound of $\E V$, and thus from a lower bound for
\begin{equation} \label{eqn:decomplower}
\E V^{\rm{idle}} = q \E V^{\rm{switching}} + (1-q) \E V^{\rm{waiting}},
\end{equation}
where $q=\pr{\rm{server~switching}|\rm{server~idle}}$.

Furthermore, we note that the cycle time satisfies
$$
(1-\rho_0)\E C = r_0 + f_0.
$$

\paragraph*{Workload while switching.} Now we express $\E V^{\rm{switching}}$ in terms of the (unknown) $f_i$. This is completely analogous to the derivation for our concrete strategy with the arguments following (\ref{eqn:wlmain2}) replacing $T_i$ by $f_i$. The result is
\begin{eqnarray}
\E V^{\rm{switching}} &= & \frac{r_{0}(\rho_{0}^{2}-\sum_{i=1}^N \rho_{i}^{2})}{2\pr{{\rm switching}}}+\frac{1}{\E C \cdot \pr{{\rm switching}}}\left\lbrace\rho_{0}\frac{r_{0}^{(2)}}{2} \right.\nonumber \\
 &  &+\left. \sum_{i=1}^N r_{i}\left[(\rho_{0}-\rho_{i})f_{i}+\sum_{j=1}^{i-1} f_{j}\left(\sum_{l=i+1}^N \rho_{l}+\sum_{l=1}^{j-1}\rho_{l}\right)+\sum_{j=i+1}^N f_{j}\sum_{l=i+1}^{j-1}\rho_{l}\right]\right\rbrace.\label{eq:EVschaltetSchranke}
\end{eqnarray}

\paragraph*{Workload while waiting.} Since we do not know the distribution of the waiting time, we cannot say much about the workload generated while the server is waiting. So, we will have to estimate at this point. First note that
\begin{equation} \label{eqn:Vlowerwaiting}
\E V^{\rm{waiting}}  =  \frac{1}{\pr{{\rm waiting}}}\,\sum_{i=1}^N p_{i}\E V_i^{\rm{waiting}},
\end{equation}
where $\E V_i^{\rm{waiting}}$ is the expected workload that is present in the system at a point in time when the server is waiting at station $i$ and $p_i=f_i/\E C$.

We cannot calculate the workload present at a point in time when we encounter the server waiting, $\E V_i^{\rm{waiting}}$, for an arbitrary strategy, but we \emph{can} give a lower bound. Namely, we can say that $\E V_i^{\rm{waiting}}$ must be at least, on the one hand, the traffic that was accumulated at the other stations during the time that the server has already passed waiting at station $i$ (that is, the expected backwards recurrence time). Additionally, since the decision of the strategy does not depend on the recent times the server has worked at the other stations nor the switchover times, we can also count the terms (\ref{eqn:wrkldworking}) and (\ref{eqn:wrkldswitching}), that is the traffic that was accumulated (and is still present) while the server was switching and working at other stations, respectively. This gives
$$
\E V_i^{\rm{waiting}} \geq \frac{f_i^{(2)}}{2 f_i} \cdot (\rho_0-\rho_i) + \sum_{k\neq i} \left( (\rho_k \E C + f_k) \cdot \sum_{j\in\{i+1,\ldots,k-1\}} \rho_j + r_k \cdot \sum_{j\in\{i+1,\ldots,k\}} \rho_j\right).
$$ 
This is the crucial observation in the derivation of the lower bound.

Now, by Jensen's inequality (\ref{eqn:jensen}), the last term can be yet bounded below by
\begin{equation} \label{eqn:vwaitli}
\frac{f_i^{2}}{2 f_i} \cdot (\rho_0-\rho_i) + \sum_{k\neq i} \left( (\rho_k \E C + f_k) \cdot \sum_{j\in\{i+1,\ldots,k-1\}} \rho_j + r_k \cdot \sum_{j\in\{i+1,\ldots,k\}} \rho_j\right).
\end{equation}
Furthermore, we need that
$$
\frac{q}{\pr{{\rm switching}}}=\frac{1-q}{\pr{{\rm waiting}}} = \frac{1}{\pr{{\rm idle}}}=\frac{1}{1-\rho_0}.
$$

Then, putting (\ref{eqn:vwaitli}) back into (\ref{eqn:Vlowerwaiting}), and this and (\ref{eq:EVschaltetSchranke}) back into (\ref{eqn:decomplower}) gives a lower bound for $\E V$ (and thus for $\bar{D}$) only in terms of the $f_i$. Minimizing over the $f_i$ leads to (\ref{eq:Wschranke2}).

Note that we cannot count the terms (\ref{eqn:terma}) and (\ref{eqn:termb}), since e.g.\ $Z_i$ and $F_i$ are not independent.

\section{Different strategies and outlook}\label{sec:outlook}
There is another strategy, which we would like to propose here (we will refer to it as ``Strategy~II''), which is likely to be better than the one proposed so far (called ``Strategy~I'' in this section) in terms of mean average delay. However, we are not able to analyse the mean average delay of Strategy~II with the present methods.

\begin{figure}
\includegraphics[scale=0.88]{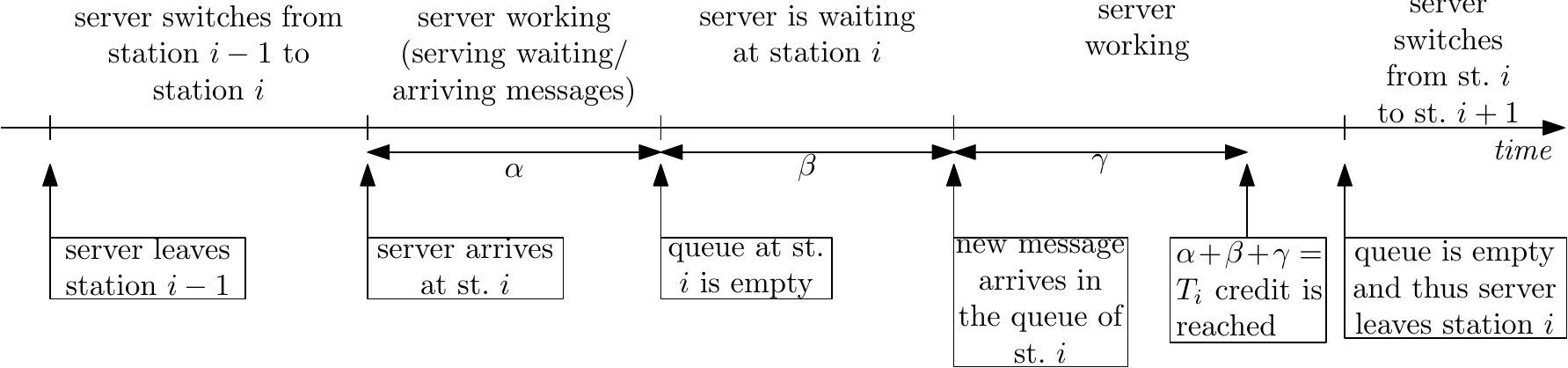} 
\caption{Operation of the polling model with Strategy II} \label{fig:stratII}
\end{figure}
Strategy~II is defined as follows. We consider a polling model as above, the only difference being that the credit $T_i$ now refers to the total time the server spends at station $i$. More precisely, a server arrives at station $i$. It then works or waits at station $i$ depending whether messages are present or not. At time $T_i$ after its arrival at the station, it only finishes all the work that is present at that time (exhaustively). It does not turn idle again at that station in the current cycle (i.e.\ into ``wait-and-see''); whenever the station is empty, it starts switching to the next station. We refer to Figure~\ref{fig:stratII} for an illustration.


It is likely that Strategy~II adjusted to its optimal waiting parameters gives a lower average delay than Strategy~I adjusted to its optimal waiting parameters $T_1^*,\ldots, T_N^*$. Heuristically, Strategy~II uses more information about the system, because it also counts the busy periods at the current station. However, we remark that even the determination of the cycle time, as in (\ref{eqn:dhash}), does not seem to be possible in a straightforward way.

We conjecture that for $N=2$ and deterministic switchover times Strategy~II provides the lowest mean average delay in the class of all strategies that are not allowed to use information of the queue status at the other station nor to look into the future of the system.

\begin{figure}
\begin{center}
\includegraphics[viewport=155 440 580 720, scale=0.9, clip=true]{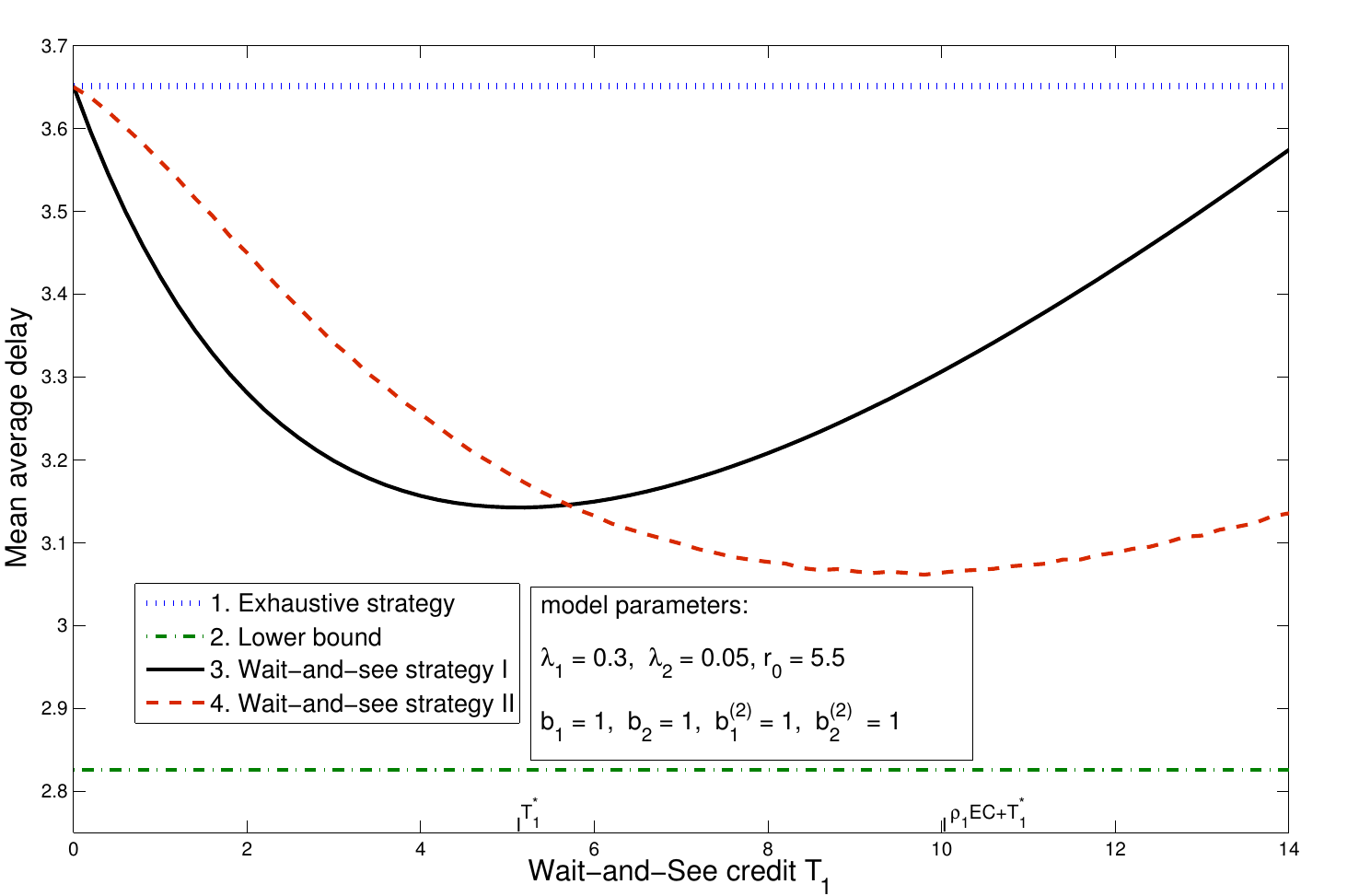}
\end{center}
\caption{Comparison of the optimal credits with Strategy~I and~II \label{fig:vglcreditIV}}
\end{figure}
Figure~\ref{fig:vglcreditIV} gives a comparison of Strategy~I and~II, where the curve for Strategy~II is obtained from simulations. We observed that the optimal credit for Strategy~II is \emph{approximated} by $T_1^*+\rho_1 \E C$, where $T_1^*$ is the optimal parameter of Strategy~I and $\E C$ is the cycle time of Strategy~I for this optimal parameter.

\medskip
Certainly, one can define different strategies, where e.g.\ the server additionally has more information on the current queue status at the other stations. This may give an average delay that is even below the lower bound given in Theorem~\ref{thm:lower}. However, note that even if the server is aware of the queue status at all stations, it is not completely clear what is the best decision at each moment in terms of lower average delay: switch or wait-and-see...

Further, one can imagine a situation where the server may look into the close future of the incoming traffic at the present station; and it may thus decide to abandon the station before the end of its wait-and-see period, when it is clear that no traffic will arrive during that time.

\paragraph*{Acknowledgement.} This work was supported by the DFG Research Center ``{\sc Matheon} -- Mathematics for key technologies'' in Berlin.

\bibliographystyle{plain}

\end{document}